# Projected likelihood contrasts for testing homogeneity in finite mixture models with nuisance parameters

Debapriya Sengupta[*,1] and Rahul Mazumder[2]

*Indian Statistical Institute*

**Abstract:** This paper develops a test for homogeneity in finite mixture models where the mixing proportions are known a priori (taken to be 0.5) and a common nuisance parameter is present. Statistical tests based on the notion of Projected Likelihood Contrasts (PLC) are considered. The PLC is a slight modification of the usual likelihood ratio statistic or the Wilk's $\Lambda$ and is similar in spirit to the Rao's score test. Theoretical investigations have been carried out to understand the large sample statistical properties of these tests. Simulation studies have been carried out to understand the behavior of the null distribution of the PLC statistic in the case of Gaussian mixtures with unknown means (common variance as nuisance parameter) and unknown variances (common mean as nuisance parameter). The results are in conformity with the theoretical results obtained. Power functions of these tests have been evaluated based on simulations from Gaussian mixtures.

## 1. Introduction

Finite mixture models are often used to understand whether the data comes from a heterogeneous or a homogeneous population. In particular, consider the case of a mixture of two populations with the mixing proportions known (Goffinet et al. [7]). We are interested to know whether the data is sampled from a proper mixture of two distributions or a single distribution.

In particular, consider a mixture family $g$, with generating population densities given by $\mathcal{M}_0 = \{f(\cdot|\theta, \eta) : \theta \in \Theta, \eta \in \mathcal{E}\}$, where $\theta$ is the main parameter of interest and $\eta$ is the common nuisance parameter. We assume that the mixing proportion is known a priori to be 0.5. The mixture model then becomes

$$(1.1) \qquad g(z|\theta_1, \theta_2, \eta) = 0.5\, f(z|\theta_1, \eta) + 0.5\, f(z|\theta_2, \eta).$$

The null hypothesis for homogeneity is, $\theta_1 = \theta_2$.

In several practical examples (for example, arising in speech analysis and non-parametric regression methodology) detection of the location of *discontinuity* in the local mean or the local variance (or local amplitude) are of interest (Figure 1). The theoretical results developed in this paper can be used in such problems. Figure 1 demonstrates several scenarios of signals being scanned through a running window of specified bandwidth. When the center of the window is placed at points of discontinuity the raw signal values ($y$-axis) will have a distribution which can be adequately modeled by (1.1). This basic idea has been explored by Hall and Titterington [8] in the context of edge and peak preserving smoothers.

---

*Supported in part by Grant No. 12(30)/04-IRSD, DIT, Govt. of India.
[1]Bayesian and Interdisciplinary Research Unit, Indian Statistical Institute, 203 B. T. Road, Kolkata 700108, India, e-mail: dps@isical.ac.in
[2]C. V. Raman Hall, 205 B. T. Road, Kolkata 700108, India, e-mail: rahul.mazumder@gmail.com
*AMS 2000 subject classifications:* Primary 62G08, 60G35; secondary 60J55.
*Keywords and phrases:* Gaussian mixture models, projected likelihood contrast.

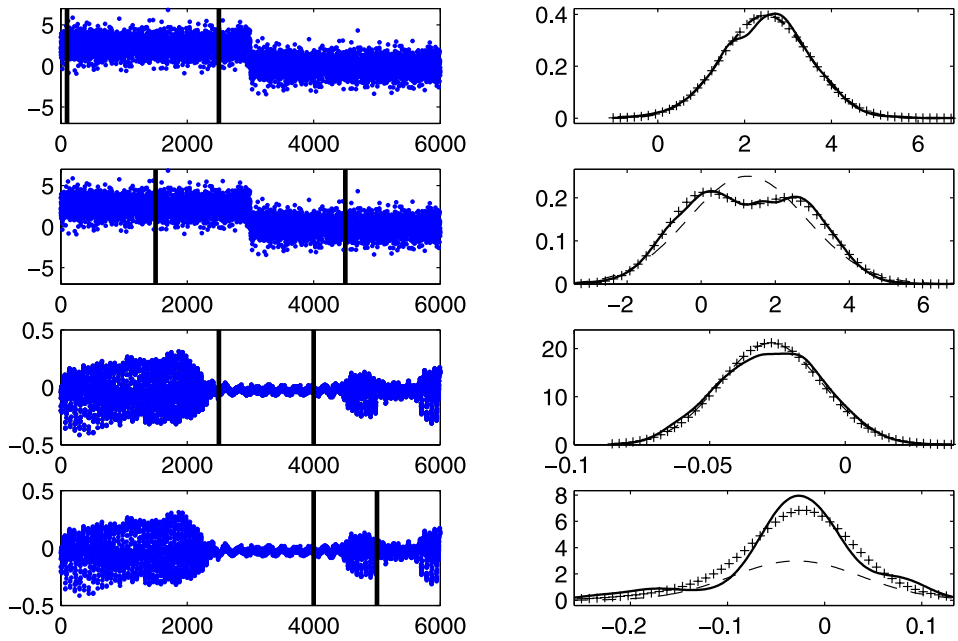

FIG 1. *Left column shows time plots of data with solid vertical lines marking the windows considered. The top two panels indicate a simulated noisy signal (with additive Gaussian noise) with mean function having a jump discontinuity. The bottom panels describe a portion of digitized speech waveform. In the right column three fitted densities of y-values: nonparametric kernel smoothed density (solid line), single component Gaussian fit (dashed line) and mixture of two Gaussian fit with equal mixing weights (curve indicated by +), are shown corresponding to the frames indicated in the left column.*

A brief list of references dealing with the study of mixture distributions and properties of the Likelihood Ratio Test (LRT) tests are provided below. In Titterington et al. [13], McLachlan and Basford [11] and Lindsay [10] one may find extensive discussions about the background of finite mixture models. The asymptotic distributions of the LRT in mixture models have been studied in Bickel and Chernoff [1], Chernoff and Lander [5], Ghosh and Sen [6], Lemdani and Pons [9]. Different modifications of LRT tests in mixture models are proposed and studied by Chen et al. [4] and Self and Liang [12].

In this paper we introduce a concept of Projected Likelihood Contrasts (PLC), a modified version of the LRT test or the Wilks' $\Lambda$ (Wilks [14]) statistic, which we motivate as follows. Consider i.i.d. observations $Z_1, Z_2, \ldots, Z_N$ generated by some element of the class of densities $g$ given by (1.1). The likelihood under the full mixture model is given by

$$(1.2) \qquad L_N(\theta_1, \theta_2, \eta) = \sum_{i=1}^{N} \log g(Z_i|\theta_1, \theta_2, \eta),$$

where $g$ is defined through (1.1). Under the null hypothesis the likelihood reduces to the usual likelihood under $\mathcal{M}_0$, namely,

$$(1.3) \qquad L_N(\theta, \theta, \eta) = \sum_{i=1}^{N} \log f(Z_i|\theta, \eta).$$

Define $(\hat{\theta}, \hat{\eta})$ as the maximum likelihood estimators of $(\theta, \eta)$ under (1.3). The idea behind the PLC statistics is to plug in the estimated nuisance parameter under the null in (1.2) and maximize it over remaining parameters $\theta_1$ and $\theta_2$. Finally the PLC statistic is defined as

$$(1.4) \qquad \Lambda_N = 2\left(\max_{\theta_1, \theta_2} L_N(\theta_1, \theta_2, \hat{\eta}) - L_N(\hat{\theta}, \hat{\theta}, \hat{\eta})\right).$$

The term projected likelihood is used here to distinguish the procedure from profile likelihood. We call it projected likelihood because the profile of the nuisance parameter is obtained after projecting the full likelihood onto $f(\cdot|\hat{\theta}, \eta) \in \mathcal{M}_0$. That way we first obtain a projected profile of $\eta$ and then maximize it so that its estimate coincides with the maximum likelihood estimate (MLE) under the null hypothesis. Note that this procedure, in spirit, is very similar to the Rao's score test.

The paper is organized as follows. In Section 2, the large sample properties of the PLC statistics is discussed. In Section 3, some simulation studies are provided. The proof of the main theorem in Section 2 is provided in the Appendix.

## 2. Large sample approximation of PLC statistic

For the purpose of theoretical investigation we shall simplify the model further assuming that the class of densities are all one dimensional. Denote the null hypothesis by

$$(2.1) \qquad H_0^* : Z_1, Z_2, \ldots, Z_N \text{ are iid } \mathcal{M}_0$$

For notational convenience we adopt the convention that the symbol $D_x^r$ indicates $r$-th partial derivative with respect to $x$, treated as a generic argument in a function. Define the following estimated scores

$$(2.2) \qquad \hat{\xi}_r(j) = \frac{D_\theta^r f(Z_j|\hat{\theta}_N, \hat{\eta}_N)}{f(Z_j|\hat{\theta}_N, \hat{\eta}_N)},$$

for $1 \leq j \leq N$ and $r \geq 1$. Analogously define the true scores $\xi_r(j) = \frac{D_\theta^r f(Z_j|\theta, \eta)}{f(Z_j|\theta, \eta)}$ at the true parameter values under $H_0^*$. One can verify that $E_{H_0^*} \xi_r(1) = 0$ for every $r \geq 1$ in case of regular parametric families.

Note that under regularity assumptions on the model the scores are well behaved and have finite moments. For the Gaussian case all moments will be finite since the joint moment generating function of any finite set of polynomials involving $\xi_r$'s exists. Define the following mixed partial derivatives of the full likelihood $L_N$.

$$(2.3) \qquad C_{ij}^N = (D_{\theta_1} + D_{\theta_2})^i (D_{\theta_1} - D_{\theta_2})^j L_N(\hat{\theta}_N, \hat{\theta}_N),$$

where $i, j$ are nonnegative integers. Moreover, let $\bar{C}_{ij}^N = N^{-1} C_{ij}^N$. Although the quantities defined in (2.3) look quite incomprehensible they can however be expressed as linear combinations of $D_{\theta_1}^l D_{\theta_2}^m L_N(\theta_1, \theta_2)$ using the Binomial expansion.

One can establish with some effort the following. $D^i_{\theta_1} D^j_{\theta_2} \log g(z|\hat{\theta}_N, \hat{\theta}_N, \hat{\eta}_N) = \sum_{\Omega}^* a(\Omega) \prod_{r=1}^{i+j} \hat{\xi}_r^{\omega_r}(z)$, where $\sum^*$ runs over all nonnegative integral partitions $\Omega = (\omega_1, \omega_2, \ldots, \omega_{p+q})$ satisfying $\sum r \omega_r = i+j$. The coefficients $a(\Omega)$ are complicated combinatorial quantities but can be recursively computed. It can be verified that $\bar{C}^N_{ij} = 0$ if $j$ is odd. We provide simplified expressions for some of the lower order $\bar{C}^N_{ij}$ which are necessary for future calculations.

(2.4)
$$\bar{C}^N_{20} = \frac{1}{N} \sum_{j=1}^N (\hat{\xi}_2(j) - \hat{\xi}_1^2(j)) \; (\xrightarrow{P} -\mathcal{I}), \quad \bar{C}^N_{02} = \frac{1}{N} \sum_{j=1}^N \hat{\xi}_2(j)$$
$$\bar{C}^N_{12} = \frac{1}{N} \sum_{j=1}^N (\hat{\xi}_3(j) - \hat{\xi}_1(j) \hat{\xi}_2(j)), \quad \bar{C}^N_{04} = \frac{1}{N} \sum_{j=1}^N \hat{\psi}(j),$$

where $\hat{\psi}(j) = \hat{\xi}_4(j) + \frac{1}{2} \hat{\xi}_1(j) \hat{\xi}_3(j) - 3 \hat{\xi}_2^2(j) + 3 \hat{\xi}_1^2(j) \hat{\xi}_2(j)$, and $\mathcal{I}$ is the Fisher information of $\theta$ under $H_0^*$. Finally, let $C_{ij}$ denote asymptotic expected values of $\bar{C}^N_{ij}$ under $H_0^*$ which can be easily derived using Lemma 2.1 (i). The distributional properties of $\bar{C}^N_{ij}$ can be derived using classical properties of $M$-estimators. We state the following lemma for the sake of completeness. The proof can be found in Bickel and Doksum [2].

**Lemma 2.1.** *Let $Z_1, Z_2, \ldots, Z_N$ be independent and identically distributed random variables with density $f(z|\theta)$ satisfying usual regularity conditions with the score function $S(z, \theta)$ and Fisher information matrix $\mathcal{I} = Cov_\theta (S(Z_1, \theta))$.*

(i) *Let $\psi(z, \theta)$ be a real valued, continuously differentiable (in $\theta$) kernel with $E_\theta \psi^2(Z_1, \theta) < \infty$, for every $\theta$. Further let $\hat{\theta}_N$ denote the MLE of $\theta$. Then*

$$\frac{1}{N} \sum_{i=1}^N \psi(Z_i, \hat{\theta}_N) \xrightarrow{P} E_\theta \psi(Z_1, \theta).$$

(ii) *In addition if $\psi$ satisfies $E_\theta \psi(Z_1, \theta) = 0$ for every $\theta$ then*

(2.5)
$$\frac{1}{\sqrt{N}} \sum_{i=1}^N \psi(Z_i, \hat{\theta}_N) \Longrightarrow N(0, V^2),$$

*where $V^2 = E_\theta \psi^2 - C' \mathcal{I}^{-1} C$ where $C = Cov_\theta (\psi(Z_1, \theta), S(Z_1, \theta))$.*

Finally, we proceed to the main asymptotic representation theorem of the PLC statistic. It turns out that even in the Gaussian case the standard $\chi^2$-approximation does not hold. Actually it turns out that Gaussian case is more paradoxical than one would expect. As a result one has to go for higher order expansion to get an idea of the limiting behavior of the statistic. The crucial issue is whether $E_{H_0^*} \xi_1(1) \xi_2(1) = 0$ or not. This is a measure of some type of spurious non-degeneracy in the model due to skewness and its asymptotic effect needs to be corrected for. Two cases are considered in the simulation section. In the first case we consider a mixture Gaussian with different means but common unknown variance and the in second case scale mixture Gaussian with common unknown mean is considered. In both cases we find $E_{H_0^*} \xi_1(1) \xi_2(1) = 0$. The first case is covered by Theorem 2.2(i) below while the second case is covered by Theorem 2.2(ii). We state the theorem keeping these two special cases in mind. The proof of the theorem is provided in the Appendix.

**Theorem 2.2.** *Assume that $E_{H_0^*} \xi_1(1) \xi_2(1) = 0$ and $C_{04} < 0$. Then under $H_0^*$,*

(i) *if $\bar{C}^N_{02} = 0$, then $\Lambda_N \xrightarrow{P} 0$.*

(ii) if $\sqrt{N}\,\bar{C}_{02}^N \Longrightarrow N(0,\sigma^2)$ *for some* $\sigma^2 > 0$, *then*

$$\Lambda_N \Longrightarrow c^2 \max(0, Z)^2, \tag{2.6}$$

*for suitable* $c^2 > 0$ *and a standard normal variate* $Z$.

## 3. Simulation studies in the case of Gaussian mixtures

In this section we provide results pertaining to the sampling distributions of the PLC statistic under the null in case of Gaussian mixtures [7]. Studies have been carried out for two different cases: unknown variances and common mean as the nuisance parameter and unknown means and common variance as the nuisance parameter. The simulation results are in conformity with the theoretical results derived. The power function of the PLC test statistic for each of the above two set-ups have been studied for different values of the alternative. Simulation studies have been carried out for different sample sizes.

### 3.1. Null distributions of the PLC

Consider the particular example of Gaussian mixture models, the main parameters of interest are the unknown means and the common variance is the nuisance parameter. The generating model is given by

$$f(z|\theta,\eta) = \eta^{-1}\phi((z-\theta)/\eta) \tag{3.1}$$

where $\phi$ is the standard normal probability density function ($\theta \in \Re, \eta > 0$). In this case $\hat{\eta}^2 = N^{-1}\sum_{i=1}^{N}(Z_i - \bar{Z})^2$, where $\bar{Z} = N^{-1}\sum_{i=1}^{N} Z_i$. The corresponding PLC is denoted by $\Lambda_N^m$. Simulation studies for the null distribution of $\Lambda_N^m$ have been performed for sample sizes $N=50$, 100 and 200. Percentiles of the sampling distribution are displayed in Table 1 which shows how different percentiles $p$ (5, 50 and 95) of the null distribution of $\Lambda_N^m$ decrease with increasing sample size $N$. The difference of the percentile values, (say that between percentiles 95 and 5), decreases with increasing sample size as well. The tabulated values give sufficient reason to believe in the validity of the theoretical results obtained in Theorem 2.2.

In the second example, also pertaining Gaussian mixture models, the main parameters of interest are unknown variances and the common mean is the nuisance parameter.

$$f(z|\theta,\eta) = \theta^{-1}\phi((z-\eta)/\theta) \tag{3.2}$$

for $\theta > 0, \eta \in \Re$. Here $\hat{\eta} = \bar{Z}$. The corresponding PLC statistic is denoted by $\Lambda_N^s$.

Table 1
*Percentiles of the null distribution of the PLC, corresponding to a Gaussian mixture with unknown means and common variance as the nuisance parameter*

|       | Percentiles |       |       |
|-------|-------|-------|-------|
| $N$   | 5     | 50    | 95    |
| 50    | 0.008 | 0.011 | 0.014 |
| 100   | 0.004 | 0.005 | 0.006 |
| 200   | 0.002 | 0.002 | 0.003 |

The constant $c^2$ in the limiting distribution (2.6) can be computed, but the computations are quite cumbersome. Hence the constant $c^2$ has been evaluated based on the sampling distribution of $\Lambda_N^s$ under the null. The sampling distribution is based on 5000 simulations of data-size 2000. The value of $c^2$ hence obtained is 0.69070.

The asymptotic null distribution of $\Lambda_N^s$ is a mixture of a degenerate mass at 0 and a $c^2 \chi_1^2$ (for suitable $c^2 > 0$), with mixing proportion 0.5. The sampling distribution of $\Lambda_N^s$, obtained from 5000 simulations of sample size 2000, is found to be a mixture of outcomes which are exactly zero and another strictly positive absolutely continuous distribution. We have observed that this absolutely continuous distribution (as obtained from simulations) is very close to $c^2 \chi_1^2$ (where $c^2 = 0.69070$) as depicted in Figure 2. Hence simulation studies of the null distribution show sufficient conformity to the theoretical results obtained in Theorem 2.2.

Simulation studies for the null distribution of $\Lambda_N^s$ have been performed and tabulated (see Table 2) for different sample sizes $N$ based on 1000 simulations of data size $N$ where $N = 50, 100, 200$.

The expected value of the sampling distribution shows a negative bias. The degree to which it approximates the mean of the large sample distribution of the PLC improves with increasing sample size. The proportion of zeros in the sampling

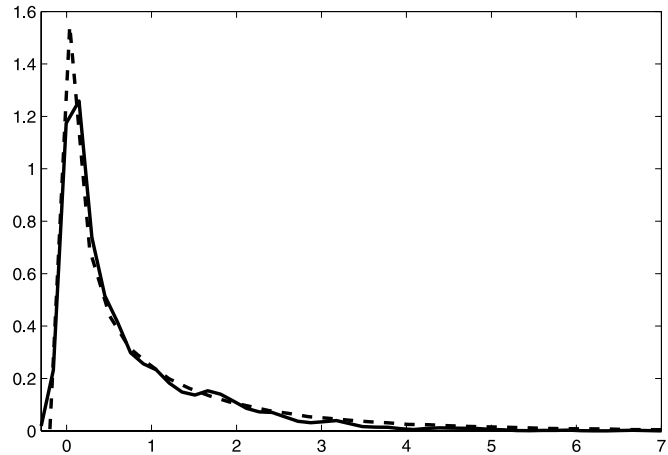

FIG 2. *Dotted line shows the kernel density estimate of $c^2(\max\{0, N(0,1)^2\})(c^2 = 0.69070)$, the theoretical asymptotic null distribution of the PLC under $N(0,1)$. Note that by invariance the results do not depend on the choice of the mean and variance. The solid line is the kernel density estimate of the sampling distribution of the PLC with the zeros left out, under the null corresponding to a Gaussian mixture of the same set-up. This sampling distribution is based on 5000 simulations of sample size 2000.*

TABLE 2
*Summary statistics of the null distribution of the PLC, corresponding to a Gaussian mixture with unknown variance and common mean as the nuisance parameter*

|     | Expectation | | % of zeros | | 5% signif. point | |
| --- | --- | --- | --- | --- | --- | --- |
| $N$ | Theor.* | Est. | Theor. | Est. | Theor.* | Est. |
| 50  | 0.345 | 0.156 | 50 | 70.1 | 1.86 | 0.935 |
| 100 | 0.345 | 0.256 | 50 | 61.5 | 1.86 | 1.608 |
| 200 | 0.345 | 0.328 | 50 | 57.5 | 1.86 | 1.817 |

*The sampling distribution based on 5000 simulations of sample-size 2000, has been used as a proxy for the theoretical asymptotic null distribution.

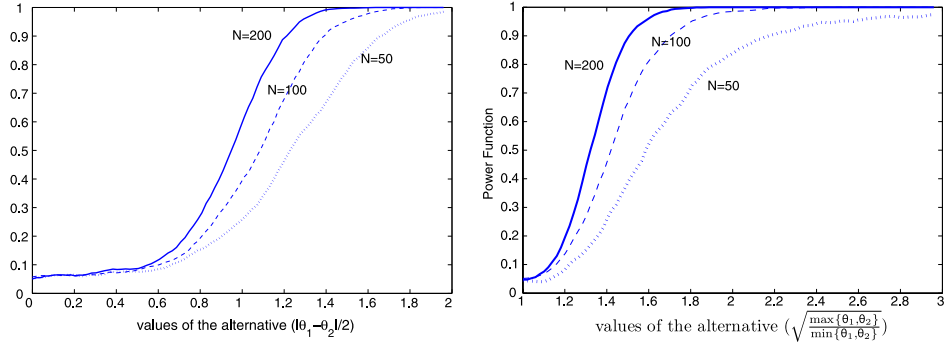

Fig 3. *Solid line, dotted line and dashed line correspond to the sample sizes 200, 100 and 50 respectively in both the figures. Power functions of the PLC test statistic at level $\alpha = 0.05$ have been evaluated. In the case of $\Lambda_N^m$, (left figure) the power function has been evaluated for values of the parameter $\frac{|\theta_1 - \theta_2|}{2} \in [0, 2]$. The power function corresponding to $\Lambda_N^s$ (right figure) has also been evaluated for the values of the parameter $\sqrt{\frac{\max\{\theta_1, \theta_2\}}{\min\{\theta_1, \theta_2\}}} \in [1, 3]$.*

distribution goes on decreasing with $N$ before it asymptotes to the theoretical value 0.5. The degree to which the sampling distribution approximates the theoretical distribution improves with increasing sample size in the case of the $95^{th}$ percentile.

### 3.2. Power function of the PLC test statistic

Power functions corresponding to the test statistic $\Lambda_N^m$ at level $\alpha = 0.05$ have been evaluated for different values of the parameter (different values of the alternative) $\frac{|\theta_1 - \theta_2|}{2}$ in the range $[0, 2]$, for three different sample sizes $N = 50, 100, 200$. (Figure 3). The power is found to increase with increasing sample size.

Power functions corresponding to the test statistic $\Lambda_N^s$ at level $\alpha = 0.05$ have been evaluated for different values of the parameter (different values of the alternative) $\sqrt{\frac{\max\{\theta_1, \theta_2\}}{\min\{\theta_1, \theta_2\}}}$ in the range $[1, 3]$, for three different sample sizes $N = 50, 100, 200$. (Figure 3). The power is found to increase with increasing sample size.

### Appendix: Proof of Theorem 2.2

First, it follows from Chen et al. [4] that both the MLEs $\hat{\theta}_1$ and $\hat{\theta}_2$ respectively are $N^{1/4}$ consistent under (1.1). For both the cases in the theorem we re-parametrize the problem with $\theta_1 = \hat{\theta}_N + N^{-1/2} s + N^{-1/4} \tau$ and $\theta_2 = \hat{\theta}_N + N^{-1/2} s - N^{-1/4} \tau$ and study its behavior near $(\hat{\theta}_N, \hat{\theta}_N)$ in the neighborhoods $|s| \leq \log N$ and $|\tau| \leq \log N$ respectively. In what follows we do not verify orders of remainder terms explicitly. Several technical steps need to be verified in the process of deriving the result. We refer to Bickel and Doksum [2], Ghosh and Sen [6] and Bose and Sengupta [3] for the type of regularity assumptions and machinery needed for uniform approximations in such a context. Also, note that under the above parametrization the likelihood becomes an even function in $\tau$. Therefore we work with $\tau \geq 0$ without any loss of generality. The asymptotic problem is non-standard because the Fisher information matrix, $\mathcal{I}(\theta_1, \theta_2, \eta)$, has rank 2 if $\theta_1 = \theta_2$ and 3 otherwise (can be verified by straightforward differentiation). Next define

$$H(s,\tau) = L_N(\hat{\theta}_N + s + \tau, \hat{\theta}_N + s - \tau).$$

It can be readily verified from (1.2) and (2.3) that

(A.1) $$\frac{\partial^{i+j}}{\partial s^i\, \partial \tau^j} H(0,0) = C^N_{ij},$$

for $i, j \geq 0$. The strategy of the proof is the following. Since the expansion is regular in *within-model displacement* $s$, we fix $\tau \geq 0$ and maximize over $s$ in the first step. Then, we examine the behavior of the maximum value obtained in the first step across $\tau$ to derive the final result. Because of our general regularity conditions all the following calculations will be valid uniformly in probability over the compact set $|s| \leq \log N$ and $0 \leq \tau \leq \log N$. In what follows $\gamma > 0$ shall denote a generic constant whose value may be determined on a case by case basis. Also, in deriving the orders of remainders we specially mention one simple fact from calculus, namely, $N^{-a}(\log N)^b \to 0$ as $N \to \infty$ for any $a, b > 0$.

(A.2) $$\begin{aligned} H(N^{-1/2}s, \tau) &= H(0,\tau) + s\,[N^{-1/2} H_{10}(0,\tau)] \\ &\quad + \frac{1}{2} s^2\,[N^{-1} H_{20}(0,\tau)] + o_P(N^{-\gamma}), \end{aligned}$$

where $H_{ij}$'s denote respective partial derivatives of $H$. Also, it can be checked that

$$H_{20}(0,\tau) = -N\mathcal{I}\,(1 + o_P(N^{-\gamma})).$$

Therefore, in large samples, for fixed $0 \leq \tau \leq N^{-1/4} \log N$, the maximum value of $H(N^{-1/2}s, \tau)$ over the compact set $|s| \leq \log N$ cannot exceed its unrestricted global maximum, which is of the order of $[N^{-1/2} H_{10}(0,\tau)]^2 / [N^{-1} H_{20}(0,\tau)]$. By direct Taylor series of order 4 we find

$$H_{10}(0, N^{-1/4}\tau) = (2!)^{-1}[\sqrt{N}\bar{C}^N_{12}]\tau^2 + (4!)^{-1}[\bar{C}^N_{14}]\tau^4 + o_P(N^{-\gamma}).$$

The facts required for the above simplification are: (i) $H_{10}(0,0) = 0$ by the maximum likelihood equation, (ii) $H_{1j}(0,0) = 0$ for $j$ odd (since $H$ is an even function of $\tau$) and (iii) the assumption of the theorem that $E_{H_0^*}\xi_1(1)\xi_2(1) = 0$. It can be checked that the last assertion implies $\sqrt{N}\bar{C}^N_{12} = O_P(1)$, in view of (2.4) and Lemma 2.1.

Therefore by virtue of the assumptions of the theorem the profile global maximum of $H(\cdot, \tau)$ becomes negligible in probability over the range of interest. Thus we have

(A.3) $$\max_{|s| \leq \log N} H(N^{-1/2}s, \tau) = H(0,\tau) + o_P(N^{-\gamma}),$$

uniformly over $0 \leq \tau \leq N^{-1/4} \log N$. Finally,

$$H(0, N^{-1/4}\tau) = H(0,0) + (2!)^{-1}[\sqrt{N}\bar{C}^N_{02}]\tau^2 + (4!)^{-1}[\bar{C}^N_{04}]\tau^4 + o_P(N^{-\gamma}).$$

Therefore we have

(A.4) $$\begin{aligned} \Lambda_N &\approx 2 \max_{|s| \leq \log N, 0 \leq \tau \leq \log N} [H(N^{-1/2}s, N^{-1/4}\tau) - H(0,0)] \\ &= \max_{0 \leq \tau \leq \log N} \{[\sqrt{N}\bar{C}^N_{02}]\tau^2 + \tfrac{1}{12}[\bar{C}^N_{04}]\tau^4 + o_P(N^{-\gamma})\}. \end{aligned}$$

Now we consider case (i) of the theorem where $\bar{C}_{02}^N = 0$. Then (A.4) reduces to $\Lambda_N = \max_{0 \leq \tau < \log N}\{(1/12)\,[\bar{C}_{04}^N]\,\tau^4 + o_P(N^{-\gamma})\}$. Since $C_{04} < 0$ it follows form Lemma 2.1 that $\Pr\{\bar{C}_{04}^N < -\delta\} \to 1$ for arbitrarily small $\delta > 0$. By choosing $\tau > 12^{1/4}\,\delta^{-1/4} N^{-\gamma/4}$ one can show that the value of the objective function (being maximized) becomes negative. Hence it can be verified that $\Lambda_N \xrightarrow{P} 0$.

For case (ii) arguing in a similar line and collecting the dominant terms from (A.3) and (A.4) and then maximizing the dominant term with respect to $\tau$ (noting that the dominant expression is a quadratic in $\tau^2$ and $\bar{C}_{04}^N \xrightarrow{P} C_{04}\,(<0)$) we obtain

$$
\begin{aligned}
\Lambda_N &\approx \max_{0 \leq \tau \leq \log N}\left\{\,[\sqrt{N}\bar{C}_{02}^N]\,\tau^2 + \tfrac{1}{12}\,[\bar{C}_{04}^N]\,\tau^4\right\} \\
&\approx -3\,\frac{[\max(0,\sqrt{N}\bar{C}_{02}^N)]^{\,2}}{\bar{C}_{04}},
\end{aligned}
\tag{A.5}
$$

with an error in approximation of the order of $o_P(N^{-\gamma})$ as before. Hence the second part of the the theorem follows from the assumptions.